\documentclass[11pt,reqno]{amsart}
\usepackage{amsmath,amssymb,graphicx}
\newtheorem{thm}{Theorem}[section]
\newtheorem{cor}[thm]{Corollary}
\newtheorem{lem}[thm]{Lemma}
\newtheorem{prop}[thm]{Proposition}
\theoremstyle{definition}
\newtheorem{defn}[thm]{Definition}
\theoremstyle{remark}
\newtheorem{rem}[thm]{Remark}
\theoremstyle{remark}
\newtheorem{ex}[thm]{Example}
\numberwithin{equation}{section}

\begin{document}

\title{A class of remarkable submartingales}%
\author{Ashkan Nikeghbali}
\address{Laboratoire de Probabilit\'es et Mod\`eles Al\'eatoires, Universit\'e Pierre et Marie Curie, et CNRS UMR 7599,
175 rue du Chevaleret F-75013 Paris, France.} \curraddr{American
Institute of Mathematics 360 Portage Ave Palo Alto, CA 94306-2244
USA} \email{ashkan@aimath.org}
\date{\today}
\begin{abstract}
In this paper, we consider the special class of positive local
submartingales $\left(X_{t}\right)$ of the form:
$X_{t}=N_{t}+A_{t}$, where the measure $\left(dA_{t}\right)$ is
carried by the set $\left\{t:\;X_{t}=0\right\}$. We show that many
examples of stochastic processes studied in the literature are in
this class and propose a unified approach based on martingale
techniques to study them. In particular, we establish some
martingale characterizations for these processes and compute
explicitly some distributions involving the pair
$\left(X_{t},A_{t}\right)$. We also associate with $X$ a solution to
the Skorokhod's stopping problem for probability measures on the
positive half-line.
\end{abstract}
\maketitle

\section{Introduction}
The study of deterministic functions or stochastic processes of the
form
\begin{equation}\label{okokko}
    Y_{t}=Z_{t}+K_{t},
\end{equation}
where $Y\geq0$ is continuous, and $K$ is increasing and continuous
with $\left(dK_{t}\right)$ carried by the set
$\left\{t:\;Y_{t}=0\right\}$, has received much attention in
probability theory: equation (\ref{okokko}) is refered to as
Skorokhod's reflection equation. It plays a key role in martingale
theory: the family of Az\'{e}ma-Yor martingales, the resolution of
Skorokhod's embedding problem, the study of Brownian local times
(see \cite{revuzyor} chapter VI for more references and examples;
see also \cite{skorokhod}, \cite{Mckean}, \cite{asterisque},
\cite{AY}, \cite{ashyordoob}). It also plays an important role in
the study of some diffusion processes (see \cite{saichotanemura,
revuzyor, asterisque}) and in the study of zeros of continuous
martingales (\cite{azemayorzero}). Well known examples of such
stochastic processes are $\left(\left|M_{t}\right|\right)$, the
absolute value of some local martingale $M$, or
$\left(S_{t}-M_{t}\right)$, where $S$ is the supremum process of
$M$. Usually, in the literature, one studies the case of the
standard Brownian Motion and its local time, then translates it into
equivalent results for the Brownian Motion and its supremum process
using L\'{e}vy's equivalence theorem, and then extends it to a wide
class of continuous local martingales, using Dubins-Schwarz
theorem.\bigskip

The aim of this paper is to provide  general framework and methods,
based on martingale techniques, to deal with a large class of
stochastic processes, which can be discontinuous, and which contain
all the previously mentioned processes. In particular, we shall
extend some well known results in the Brownian setting (and whose
proofs are based on excursion theory) to a wide class of stochastic
processes which are not Markov and which can even be discontinuous.
\bigskip

More precisely, we shall consider the following class of local
submartingales, whose definition goes back to Yor \cite{yorinegal}:
\begin{defn}\label{martreflechies}
Let $\left(X_{t}\right)$ be a positive local submartingale, which
decomposes as:
$$X_{t}=N_{t}+A_{t}.$$
We say that $\left(X_{t}\right)$ is of class $(\Sigma)$ if:
\begin{enumerate}
\item $\left(N_{t}\right)$ is a c\`{a}dl\`{a}g local martingale, with $N_{0}=0$;
\item $\left(A_{t}\right)$ is a continuous increasing process, with $A_{0}=0$;
\item the measure $\left(dA_{t}\right)$ is carried by the set
$\left\{t:\;X_{t}=0\right\}$.
\end{enumerate}If additionally, $\left(X_{t}\right)$ is of class
$(D)$, we shall say that $\left(X_{t}\right)$ is of class $(\Sigma
D)$.
\end{defn} In Section 2, we prove a martingale characterization
for the processes of class $(\Sigma)$ and then give some examples.
In particular, we obtain a family of martingales reminiscent of the
family of Az\'{e}ma-Yor martingales.

\noindent In Section 3, we prove the following estimate for a large
class of processes of class $(\Sigma)$, generalizing a well known
result for the pair $\left(B_{t},\ell_{t}\right)$, where
$\left(B_{t}\right)$ is the standard one dimensional Brownian Motion
and $\left(\ell_{t}\right)$ its local time at $0$:
\begin{equation}\label{estimate1}
\mathbb{P}\left(\exists t\geq
0,\;X_{t}>\varphi\left(A_{t}\right)\right)=1-\exp\left(-\int_{0}^{\infty}\frac{dx}{\varphi\left(x\right)}\right),
\end{equation}where $\varphi$ is a positive Borel function. We then use this estimate to obtain
the law of the maximum of some processes involving the pair
$\left(X_{t},\varphi\left(A_{t}\right)\right)$. We also use the
domination relations of Lenglart \cite{lenglart} to obtain some
probabilistic inequalities involving the pair
$\left(X_{t},\varphi\left(A_{t}\right)\right)$.

\noindent In Section 4, we compute the distribution of $A_{\infty}$
when $\left(X_{t}\right)$ is of class $(\Sigma D)$, and then deduce
the law of $A_{T}$, where $T$ is a stopping time chosen in a family
of stopping times reminiscent of the stopping times used by
Az\'{e}ma and Yor for the resolution of the Skorokhod stopping
problem. Among other applications of these results, we recover a
result of Lehoczky (\cite{lehoczky}) about the law of the maximum of
some stopped diffusions. We also give the law of the increasing
process of a conveniently stopped process which solves a stochastic
differential equation of Skorokhod type with a reflecting condition
at $0$  (such processes are defined in \cite{saichotanemura} and
play a key role in extensions of Pitman's theorem).

\noindent Finally in Section 5, inspired by a recent paper of
Obl\'{o}j and Yor (\cite{jan}), and using the estimates
(\ref{estimate1}), we propose two different solutions to the
Skorokhod embedding problem for a non atomic probability measure on
$\mathbb{R}_{+}$.

\section{A first characterization and some examples}
Let $\left(X_{t}\right)$ be of class $(\Sigma)$. We have the
following martingale characterization for the processes of class
$(\Sigma)$:
\begin{thm}\label{caracteriation1}
The following are equivalent:
\begin{enumerate}
\item The local submartingale $\left(X_{t}\right)$ is of class $(\Sigma)$;

\item There exists an increasing, adapted and continuous process $\left(
C_{t}\right) $ such that for every locally bounded Borel function
$f$, and $F\left(x\right)\equiv\int_{0}^{x}f\left(z\right)dz$, the
process
\begin{equation*}
F\left( C_{t}\right) -f\left( C_{t}\right) X_{t}
\end{equation*}%
is a local martingale. Moreover, in this case, $\left( C_{t}\right)$
is equal to $\left( A_{t}\right)$, the increasing process of $X$.
\end{enumerate}
\end{thm}
\begin{proof}
$\left( 1\right) \Longrightarrow \left( 2\right) .$ First, let us
assume that $f$ is $\mathcal{C}^{1}$ and let us take $C_{t}\equiv
A_{t}$. An integration by parts and the fact that
$\int_{0}^{t}X_{s-}dA_{s}=\int_{0}^{t}X_{s}dA_{s}$ since $A$ is
continuous yield:
\begin{eqnarray*}
  f\left( A_{t}\right) X_{t} &=& \int_{0}^{t}f\left( A_{u}\right)dX_{u}+\int_{0}^{t}f'\left( A_{u}\right)X_{u}dA_{u} \\
   &=&  \int_{0}^{t}f\left( A_{u}\right)dN_{u}+\int_{0}^{t}f\left( A_{u}\right)dA_{u}+\int_{0}^{t}f'\left(
   A_{u}\right)X_{u}dA_{u}.
\end{eqnarray*}
Since $\left(dA_{t}\right)$ is carried by the set
$\left\{t:\;X_{t}=0\right\}$, we have $\int_{0}^{t}f'\left(
   A_{u}\right)X_{u}dA_{u}=0$. As $\int_{0}^{t}f\left( A_{u}\right)dA_{u}=F\left(
   A_{t}\right)$, we have thus obtained that:
   \begin{equation}\label{representatiointegrale}
   F\left( A_{t}\right) -f\left( A_{t}\right) X_{t}=-\int_{0}^{t}f\left(
   A_{u}\right)dN_{u},
   \end{equation}and consequently $\left(F\left( A_{t}\right) -f\left( A_{t}\right)
   X_{t}\right)$ is a local martingale. The general case when $f$ is
   only assumed to be locally bounded follows from a monotone class
   argument and the integral representation (\ref{representatiointegrale}) is still valid.

$\left( 2\right) \Longrightarrow \left( 1\right) .$ First take
$F\left( a\right) =a$; we then obtain that $C_{t}-X_{t}$ is a local
martingale. Hence
the increasing process of $X$ in its Doob-Meyer decomposition is $C$, and $%
C=A$. Next, we take: $F\left( a\right) =a^{2}$ and we get:%
\begin{equation*}
A_{t}^{2}-2A_{t}X_{t}
\end{equation*}%
is a local martingale. But%
\begin{equation*}
A_{t}^{2}-2A_{t}X_{t}=2\int_{0}^{t}A_{s}\left( dA_{s}-dX_{s}\right)
-2\int_{0}^{t}X_{s}dA_{s}.
\end{equation*}%
Hence, we must have:%
\begin{equation*}
\int_{0}^{t}X_{s}dA_{s}=0.
\end{equation*}%
Thus $dA_{s}$ is carried by the set of zeros of $X$.
\end{proof}
\begin{rem}
We have proved that for $f$ a locally bounded Borel function, we
have:
$$ f\left( A_{t}\right) X_{t}-F\left( A_{t}\right)=\int_{0}^{t}f\left(
   A_{u}\right)dN_{u}.$$ When $f$ is nonnegative, this means that
   $$f\left( A_{t}\right) X_{t}=\int_{0}^{t}f\left(
   A_{u}\right)dN_{u}+F\left( A_{t}\right),$$is again of class
   $(\Sigma)$, and its increasing process is $F\left( A_{t}\right)$.
\end{rem}
\begin{rem}\label{reprsedesaydsremark}
When $X$ is continuous, we can apply Skorokhod's reflection lemma to
obtain:
\begin{equation*} A_{t}=\sup_{s\leq t}\left( -N_{s}\right).
\end{equation*}
Then an application of the balayage formula (see \cite{revuzyor},
chapter VI, p.262) yields:
\begin{equation*}\label{mlocbalay}
F\left( C_{t}\right) -f\left( C_{t}\right) X_{t}
\end{equation*}%
is a local martingale. In the continuous case, we can also note that
$X$ may be represented as $\left(S_{t}-M_{t}\right)$, where $M$ is a
continuous local martingale and $S$ its supremum process.
\end{rem}
One often needs to know when $\left(F\left( A_{t}\right) -f\left(
A_{t}\right)X_{t}\right)$ is a true martingale. We shall deal with a
special case which will be of interest to us later.
\begin{cor}\label{corolfbornee}
Let $X$ be of class $(\Sigma D)$. If $f$ is a Borel bounded function
with compact support, then $\left(F\left( A_{t}\right) -f\left(
A_{t}\right)X_{t}\right)$ is a uniformly integrable martingale.
\end{cor}
\begin{proof}
There exist two constants $C>0,\;K>0$ such that $\forall
x\geq0,\;|f\left(x\right)|\leq C$, and $\forall x \geq
K,\;f\left(x\right)=0$. Consequently, we have:
$$|F\left( A_{t}\right) -f\left(
A_{t}\right)X_{t}|\leq CK+CX_{t};$$now, since
$\left(CK+CX_{t}\right)$ is of class $(D)$, we deduce that
$\left(F\left( A_{t}\right) -f\left( A_{t}\right)X_{t}\right)$ is a
local martingale of class $(D)$ and hence it is a uniformly
integrable martingale.
\end{proof}\bigskip

Now, we shall give some examples of processes in the class
$(\Sigma)$.
\begin{enumerate}
\item Let $\left(M_{t}\right)$ be a continuous local martingale with
respect to some filtration $\left(\mathcal{F}_{t}\right)$, starting
from $0$; then,
$$|M_{t}|=\int_{0}^{t}sgn\left(M_{u}\right)dM_{u}+L_{t}\left(M\right),$$ is of
class $(\Sigma)$.
\item Similarly, for any $\alpha>0,\;\beta>0$, the process:
$$\alpha M_{t}^{+}+\beta M_{t}^{-}$$ is of
class $(\Sigma)$.
\item Let $\left(M_{t}\right)$ be a local martingale (starting from $0$) with only
negative jumps and let $S_{t}\equiv\sup_{u\leq t}M_{u}$; then
$$X_{t}\equiv S_{t}-M_{t}$$is of
class $(\Sigma)$. In this case, $X$ has only positive jumps.
\item Let $\left(R_{t}\right)$ be a Bessel process, starting from $0$, of dimension
$\delta$, with $\delta\in(0,2)$. We borrow our results about these
processes from \cite{Ashkanbessel} where all the proofs can be
found. We shall introduce the parameter $\mu\in(0,1)$, defined by
$\delta=2\left(1-\mu\right)$. Then,
$$X_{t}\equiv R_{t}^{2\mu}$$ is of class $(\Sigma)$ and can be
decomposed as:
$$R_{t}=N_{t}+L_{t}\left(R\right),$$where
$\left(L_{t}\left(R\right)\right)$ is one normalization for the
local time at $0$ of $R$.
\item Again, let $\left(R_{t}\right)$ be  a Bessel process (starting from $0$) of dimension
$2(1-\mu)$, with $\mu\in(0,1)$. Define:
$$g_{\mu}\left(t\right)\equiv\sup\left\{u\leq t:\;R_{u}=0\right\}.$$
In the filtration
$\mathcal{G}_{t}\equiv\mathcal{F}_{g_{\mu}\left(t\right)}$ of the
zeros of the Bessel process $R$, the stochastic process:
$$X_{t}\equiv \left(t-g_{\mu}\left(t\right)\right)^{\mu},$$is a submartingale of
class $(\Sigma)$, whose increasing process in its Doob-Meyer
decomposition is given by:
$$A_{t}\equiv
\dfrac{1}{2^{\mu}\Gamma\left(1+\mu\right)}L_{t}\left(R\right),$$where
as usual $\Gamma$ stands for Euler's gamma function. Recall that
$\mu\equiv\dfrac{1}{2}$ corresponds to the absolute value of the
standard Brownian Motion; thus for $\mu\equiv\dfrac{1}{2}$ the above
result leads to nothing but the celebrated second Az\'{e}ma's
martingale ($X_{t}-A_{t}$, see \cite{azemayormartrem,zurich}). In
this example, $X$ has only negative jumps.
\end{enumerate}
\section{Some estimates and distributions for the pair $\left(X_{t},A_{t}\right)$}
In the next sections, we shall introduce the family of stopping
times $T_{\varphi}$, associated with a nonnegative Borel  function
$\varphi$, and defined by:
$$T_{\varphi}\equiv\inf\left\{t\geq0:\;\varphi\left(A_{t}\right)X_{t}\geq1\right\}.$$These
stopping times (associated with a suitable $\varphi$) play an
important role in the resolution by Az\'{e}ma and Yor (\cite{AY})
and Obl\'{o}j and Yor (\cite{jan}) of the Skorokhod embedding
problem for the Brownian Motion and the age of Brownian excursions.
Some special cases of this family are also studied in
\cite{jeulinyor81} and \cite{lehoczky}. One natural and important
question is whether the stopping time $T_{\varphi}$ is almost surely
finite or not. The next theorem, reminiscent of some studies by
Knight (\cite{Knight73}, \cite{knight78}), answers this question,
and generalizes a result of Yor (\cite{columbia}) for the standard
Brownian Motion. In particular, martingale techniques will allow us
to get rid of the Markov property (which was one of the aims of Paul
Andr\'{e} Meyer when he developed systematically the general theory
of stochastic processes). Before stating and proving our main
Theorem, we shall need an elementary, yet powerful lemma, which we
have called Doob's maximal identity in \cite{ashyordoob}. For sake
of completeness, we give again a short proof for it.
\begin{lem}[Doob's maximal
identity] \label{maxeq} Let $\left(M_{t}\right)$ be a positive local
martingale which satisfies:
$$M_{0}=x,\;x>0;\;\lim_{t\rightarrow\infty}M_{t}=0.$$
If we note $$S_{t}\equiv\sup_{u\leq t}M_{u},$$and if $S$ is
continuous, then for any $a>0$, we have:
\begin{enumerate}
\item
\begin{equation} \mathbf{P}\left( S_{\infty }>a\right) =\left(
\frac{x}{a}\right) \wedge 1. \label{loimax}
\end{equation}Hence, $\dfrac{x}{S_{\infty }}$ is a uniform random variable on $%
\left(0,1\right)$.
\item For any stopping time $T$:%
\begin{equation}
\mathbf{P}\left( S^{T}>a\mid \mathcal{F}_{T}\right) =\left( \frac{M_{T}}{a}%
\right) \wedge 1 ,  \label{loimaxcond}
\end{equation}%
where
\begin{equation*}
S^{T}=\sup_{u\geq T}M_{u}.
\end{equation*}%
Hence, $\dfrac{M_{T}}{S^{T}}$ is also a uniform random variable on
$\left(0,1\right)$, independent of $\mathcal{F}_{T}$.
\end{enumerate}
\end{lem}

\begin{proof}
Formula (\ref{loimaxcond}) is a consequence of (\ref{loimax}) when
applied
to the martingale $\left( M_{T+u}\right) _{u\geq 0}$ and the filtration $%
\left( \mathcal{F}_{T+u}\right) _{u\geq 0}$. Formula (\ref{loimax})
itself is obvious when $a\leq x$, and for $a>x$, it is obtained by
applying Doob's optional stopping theorem to the local martingale
$\left( M_{t\wedge T_{a}}\right) $, where $T_{a}=\inf \left\{ u\geq
0:\text{ }M_{u}\geq a\right\} $.
\end{proof}
Now, we state the main result of this section:
\begin{thm}\label{estimationavecA}
Let $X$ be a local submartingale of the class $(\Sigma)$, with only
negative jumps, such that $A_{\infty}=\infty$. Define
$\left(\tau_{u}\right)$ the right continuous inverse of $A$:
$$\tau_{u}\equiv\inf\left\{t:\;A_{t}>u\right\}.$$ Let
$\varphi:\mathbb{R}_{+}\rightarrow\mathbb{R}_{+}$ be a Borel
function. Then, we have the following estimates:
\begin{equation}\label{estimationavecphi}
\mathbb{P}\left(\exists t\geq
0,\;X_{t}>\varphi\left(A_{t}\right)\right)=1-\exp\left(-\int_{0}^{\infty}\frac{dx}{\varphi\left(x\right)}\right),
\end{equation}and
\begin{equation}\label{estimationavecphiettpslocal}
\mathbb{P}\left(\exists t\leq
\tau_{u},\;X_{t}>\varphi\left(A_{t}\right)\right)=1-\exp\left(-\int_{0}^{u}\frac{dx}{\varphi\left(x\right)}\right).
\end{equation}
\end{thm}
\begin{proof}
The proof is based on Theorem \ref{caracteriation1} and Lemma
\ref{maxeq}. We shall first prove equation
(\ref{estimationavecphi}), and for this, we first note that we can
always assume that $\dfrac{1}{\varphi}$ is bounded and integrable.
Indeed, let us consider the event
$$\Delta_{\varphi}\equiv\left\{\exists t\geq0,\;X_{t}>\varphi\left(A_{t}\right)\right\}.$$
Now, if $\left(\varphi_{n}\right)_{n\geq1}$ is a decreasing sequence
of functions with limit $\varphi$, then the events
$\left(\Delta_{\varphi_{n}}\right)$ are increasing, and
$\bigcup_{n}\Delta_{\varphi_{n}}=\Delta_{\varphi}$. Hence, by
approximating $\varphi$ from above, we can always assume that
$\dfrac{1}{\varphi}$ is bounded and integrable.

Now, let $$F\left(x\right)\equiv
1-\exp\left(-\int_{x}^{\infty}\dfrac{dz}{\varphi\left(z\right)}\right);$$its
Lebesgue derivative $f$ is given by:
$$f\left(x\right)=\dfrac{-1}{\varphi\left(x\right)}\exp\left(-\int_{x}^{\infty}\dfrac{dz}{\varphi\left(z\right)}\right)=\dfrac{-1}{\varphi\left(x\right)}\left(1-F\left(x\right)\right).$$
Now, from Theorem \ref{caracteriation1}, $\left(M_{t}\equiv F\left(
A_{t}\right) -f\left( A_{t}\right) X_{t}\right)$, which is also
equal to $F\left( A_{t}\right) +
X_{t}\dfrac{1}{\varphi\left(A_{t}\right)}\left(1-F\left(A_{t}\right)\right)$,
is a positive local martingale (whose supremum is continuous since
$\left(M_{t}\right)$ has only negative jumps), with
$M_{0}=1-\exp\left(-\int_{0}^{\infty}\frac{dx}{\varphi\left(x\right)}\right)$.
Moreover, as $\left(M_{t}\right)$ is a positive local martingale, it
converges  almost surely as $t\rightarrow\infty$. Let us now
consider $M_{\tau_{u}}$:
$$M_{\tau_{u}}=F\left(
u\right)-f\left(u\right) X_{\tau_{u}}.$$But since
$\left(dA_{t}\right)$ is carried by the zeros of $X$ and since
$\tau_{u}$ corresponds to an increase time of $A$, we have
$X_{\tau_{u}}=0$. Consequently,
$$\lim_{u\rightarrow\infty}M_{\tau_{u}}=\lim_{u\rightarrow\infty}F\left(
u\right)=0,$$ and hence $$\lim_{u\rightarrow\infty}M_{u}=0.$$ Now
let us note that if for a given $t_{0}<\infty$, we have
$X_{t_{0}}>\varphi\left(A_{t_{0}}\right)$, then we must have:
$$M_{t_{0}}>F\left(
A_{t_{0}}\right) -f\left( A_{t_{0}}\right)
\varphi\left(A_{t_{0}}\right)=1,$$and hence we easily deduce that:
\begin{eqnarray*}
  \mathbb{P}\left(\exists t\geq
0,\;X_{t}>\varphi\left(A_{t}\right)\right) &=& \mathbb{P}\left(\sup_{t\geq0}M_{t}>1\right) \\
   &=&  \mathbb{P}\left(\sup_{t\geq0}\dfrac{M_{t}}{M_{0}}>\dfrac{1}{M_{0}}\right)\\
   &=& M_{0},
\end{eqnarray*}where the last equality is obtained by an application
of Doob's maximal identity (Lemma \ref{maxeq}).

To obtain the second identity of the Theorem, it suffices to replace
$\varphi$ by the function $\varphi_{u}$ defined as:
$$\varphi_{u}\left(x\right)=\begin{cases}
\varphi\left(x\right)&\mathrm{if\ }x<u\\
\infty&\mathrm{otherwise}. \end{cases}$$
\end{proof}
Now, as an application of Theorem \ref{estimationavecA}, we have the
following corollaries:
\begin{cor}
Let $X$ be a local submartingale of the class $(\Sigma)$, with only
negative jumps, such that $\lim_{t\rightarrow\infty}A_{t}=\infty$.
If $\int_{0}^{\infty}\frac{dx}{\varphi\left(x\right)}=\infty$, then
the stopping time $T_{\varphi}$ is finite almost surely, i.e.
$T_{\varphi}<\infty$. Furthermore, if $T<\infty$ and if $\varphi$ is
locally bounded, then:
\begin{equation}\label{phiategalx}
    X_{T}=\dfrac{1}{\varphi\left(A_{T}\right)}.
\end{equation}
\end{cor}
\begin{proof}
At this stage, only formula (\ref{phiategalx}) is not trivial. We
note, from remark (\ref{reprsedesaydsremark}) that:
$$\varphi\left(A_{t}\right)X_{t}=\int_{0}^{t}\varphi\left(A_{u}\right)dN_{u}+\Phi\left(A_{t}\right),$$where
$\Phi\left(x\right)=\int_{0}^{x}dz\varphi\left(z\right)$. Hence,
$\left(\varphi\left(A_{t}\right)X_{t}\right)$ has also only negative
jumps and consequently $\varphi\left(A_{T}\right)X_{T}=1$.
\end{proof}
\begin{rem}
It is remarkable that no more regularity than  local boundedness is
required for $\varphi$ to have
$X_{T}=\dfrac{1}{\varphi\left(A_{T}\right)}$. Usually, in the
literature, one requires left or right continuity for $\varphi$.
\end{rem}
\begin{rem}
Sometimes, in the literature (see section 5), one is interested in
the stopping time:
$$T\equiv\inf\left\{t\geq0:\;X_{t}\geq\varphi\left(A_{t}\right)\right\}.$$In
this case,  if $\int_{0}^{\infty}dx\varphi\left(x\right)=\infty$,
then $T<\infty$. In the case $T<\infty$ and  $1/\varphi$ is locally
bounded, we have:
$$X_{T}=\varphi\left(A_{T}\right).$$
\end{rem}
\begin{cor}[Knight \cite{knight78}]
Let $\left(B_{t}\right)$ denote a standard Brownian Motion, and $S$
its supremum process. Then, for $\varphi$ a nonnegative Borel
function, we have:
\begin{equation*}
  \mathbb{P}\left(\forall t\geq
0,\;S_{t}-B_{t}\leq\varphi\left(S_{t}\right)\right)=
\exp\left(-\int_{0}^{\infty}\frac{dx}{\varphi\left(x\right)}\right).
\end{equation*}Furthermore, if we let $T_{x}$ denote the stopping
time:
$$T_{x}=\inf\left\{t\geq0:\;S_{t}>x\right\}=\inf\left\{t\geq0:\;B_{t}>x\right\},$$then
for any nonnegative Borel function $\varphi$, we have:
\begin{equation*}
  \mathbb{P}\left(\forall t\leq
T_{x},\;S_{t}-B_{t}\leq\varphi\left(S_{t}\right)\right)=
\exp\left(-\int_{0}^{x}\frac{dx}{\varphi\left(x\right)}\right).
\end{equation*}
\end{cor}
\begin{proof}
It is a consequence of Theorem \ref{estimationavecA}, with
$X_{t}=S_{t}-B_{t}$, and $A_{t}=S_{t}$.
\end{proof}
\begin{rem}
The same result holds for any continuous local martingale such that
$S_{\infty}=\infty$.
\end{rem}
\begin{cor}[\cite{Ashkanbessel}]
Let $R$ be a Bessel process of dimension $2(1-\mu)$, with
$\mu\in(0,1)$, and $L$ its local time at $0$ (as defined in Section
2). Define $\tau$ the right continuous inverse of $L$:
$$\tau_{u}=\inf\left\{t\geq0;\;L_{t}> u\right\}.$$
Then, for any positive Borel function $\varphi$, we have:
$$\mathbb{P}\left(\exists t\leq
\tau_{u},\;R_{t}>\varphi\left(L_{t}\right)\right)=1-\exp\left(-\int_{0}^{u}\frac{dx}{\varphi^{2\mu}\left(x\right)}\right),$$and
$$\mathbb{P}\left(\exists t\geq
0,\;R_{t}>\varphi\left(L_{t}\right)\right)=1-\exp\left(-\int_{0}^{\infty}\frac{dx}{\varphi^{2\mu}\left(x\right)}\right).$$When
$\mu=\dfrac{1}{2}$, we recover the well known estimates for the
standard Brownian Motion and its local time $\left(\ell_{t}\right)$
(see \cite{columbia}).
\end{cor}
\begin{proof}
This follows from the fact that $X_{t}=R_{t}^{2\mu}$ is a
submartingale satisfying the condition of Theorem
\ref{estimationavecA} (see Section 2).
\end{proof}
\begin{cor}
Let $R$ be a Bessel process of dimension $2(1-\mu)$, with
$\mu\in(0,1)$, Define:
$$g_{\mu}\left(t\right)\equiv\sup\left\{u\leq t:\;R_{u}=0\right\}.$$
In the filtration
$\mathcal{G}_{t}\equiv\mathcal{F}_{g_{\mu}\left(t\right)}$,
$$X_{t}\equiv \left(t-g_{\mu}\left(t\right)\right)^{\mu},$$is a submartingale of
class $(\Sigma)$ whose increasing process is
$$\dfrac{1}{2^{\mu}\Gamma\left(1+\mu\right)}L_{t}\left(R\right).$$Consequently,
we have:
$$\mathbb{P}\left(\exists t\geq
0,\;g_{\mu}\left(t\right)<t-\varphi\left(L_{t}\right)\right)=1-\exp\left(-\frac{1}{2^{\mu}\Gamma\left(1+\mu\right)}\int_{0}^{\infty}\dfrac{dx}{\varphi^{2\mu}\left(x\right)}\right),$$
and
$$\mathbb{P}\left(\exists t\leq
\tau_{u},\;g_{\mu}\left(t\right)<t-\varphi\left(L_{t}\right)\right)=1-\exp\left(-\frac{1}{2^{\mu}\Gamma\left(1+\mu\right)}\int_{0}^{2^{\mu}\Gamma\left(1+\mu\right)u}\frac{dx}{\varphi^{2\mu}\left(x\right)}\right).$$
\end{cor}\bigskip

Now, we shall use Theorem \ref{estimationavecA} to obtain the
distribution of the supremum of some random variables involving the
pair $\left(X_{t},A_{t}\right)$. The results we shall prove have
been obtained by Yor (\cite{columbia}) in the Brownian setting.
\begin{cor}
Consider for $q>p>0$, the random variable:
$$S_{p,q}\equiv
\sup_{t\geq0}\left(X_{t}^{p}-A_{t}^{q}\right).$$Then, under the
assumptions of Theorem \ref{estimationavecA} for $X$,
$$S_{p,q}\ \stackrel{\mbox{\small (law)}}{=}\ \
\left(\mathbf{e}_{p,q}\right)^{\frac{pq}{p-q}},$$where
$\mathbf{e}_{p,q}$ is an exponential random variable with parameter
$c_{p,q}=\frac{1}{q}\int_{0}^{\infty}\frac{dz}{z^{1-\frac{1}{q}}\left(1+z\right)^{\frac{1}{p}}}.$
\end{cor}
\begin{proof}
Let $a>0$;
$$\mathbb{P}\left(S_{p,q}>a\right)=\mathbb{P}\left(\exists
t\geq
0,\;X_{t}>\left(a+L_{t}^{q}\right)^{\frac{1}{p}}\right)=1-\exp\left(-\int_{0}^{\infty}\dfrac{dx}{\left(a+x^{q}\right)^{\frac{1}{p}}}\right),$$and
the result follows from a straightforward change of variables in the
last integral.
\end{proof}
\begin{cor}
Let $\varphi$ be a nonnegative integrable function and define the
random variable $S_{\varphi}$ as:
$$S_{\varphi}=\sup_{t\geq0}\left(X_{t}-\varphi\left(A_{t}\right)\right).$$Then,
for all $a\geq0$, under the assumptions of Theorem
\ref{estimationavecA} for $X$, we have:
$$\mathbb{P}\left(S_{\varphi}>a\right)=1-\exp\left(-\int_{0}^{\infty}\dfrac{dx}{a+\varphi\left(x\right)}\right).$$
\end{cor}
\begin{proof}
It suffices to remark that
$$\mathbb{P}\left(S_{\varphi}>a\right)=\mathbb{P}\left(\exists t\geq
0,\;X_{t}>a+\varphi\left(A_{t}\right)\right),$$and then apply
Theorem \ref{estimationavecA}.
\end{proof}\bigskip

To conclude this section, we shall give some inequalities for the
pair $\left(X_{t},A_{t}\right)$. Inequalities for submartingales
have been studied in depth by Yor in \cite{yorinegal}: here we first
simply apply the domination inequalities of Lenglart
(\cite{lenglart}) to derive a simple inequality between
$\left(X_{t}^{*}\equiv\sup_{s\leq t}X_{s}\right)$ and
$\left(A_{t}\right)$ and then we mention a result of Yor
\cite{yorinegal} for $H^{p}$ norms.
\begin{defn}[Domination relation]
A positive adapted right continuous process $X$ ($X_{0}=0$) is
dominated by an increasing process $A$ ($A_{0}=0$) if
$$\mathbb{E}\left[X_{T}\right]\leq\mathbb{E}\left[A_{T}\right]$$for
any bounded stopping time $T$.
\end{defn}
\begin{lem}[\cite{lenglart}, p.173]
If $X$ is dominated by $A$ and $A$ is continuous, then for any
$k\in(0,1)$,
$$\mathbb{E}\left[\left(X_{\infty}^{*}\right)^{k}\right]\leq\dfrac{2-k}{1-k}\mathbb{E}\left[A_{\infty}^{k}\right].$$
\end{lem}Now, applying this to positive submartingales yields:
\begin{prop}
Let $X$ be a positive local submartingale which decomposes as:
$X_{t}=N_{t}+A_{t}$. Then, we have for any $k\in(0,1)$:
$$\mathbb{E}\left[\left(X_{\infty}^{*}\right)^{k}\right]\leq\dfrac{2-k}{1-k}\mathbb{E}\left[A_{\infty}^{k}\right].$$
If furthermore the process $\left(X_{t}^{*}\right)$ is continuous
(this is the case if $X$ has only negative jumps), then we have for
any $k\in(0,1)$:
$$\mathbb{E}\left[\left(X_{\infty}^{*}\right)^{k}\right]\leq\dfrac{2-k}{1-k}\mathbb{E}\left[A_{\infty}^{k}\right]\leq\left(\dfrac{2-k}{1-k}\right)^{2}\mathbb{E}\left[\left(X_{\infty}^{*}\right)^{k}\right].$$
\end{prop}
\begin{proof}
The proof is a consequence of:
$$\mathbb{E}\left[X_{T}\right]\leq\mathbb{E}\left[A_{T}\right]\leq\mathbb{E}\left[X_{T}^{*}\right]$$for
any bounded stopping time, combined with an application of
Lenglart's domination results.
\end{proof}
\begin{rem}
The previous proposition applies in particular to the pair
$\left(|B_{t}|, \ell_{t}\right)$, where $B$ is the standard Brownian
Motion and $\left(\ell_{t}\right)$ its local time at $0$. More
generally, it holds for the pair
$\left(R_{t}^{2\mu},L_{t}\left(R\right)\right)$ where $R$ is a
Bessel process of dimension $2(1-\mu)$. It also applies to the age
process and its local time at zero:
$\left(\left(t-g_{\mu}\left(t\right)\right)^{\mu},\frac{1}{2^{\mu}\Gamma\left(1+\mu\right)}L_{t}\left(R\right)\right)$.
More precisely, for any stopping time $T$, and any $k\in(0,1)$, we
have:
\begin{eqnarray*}
  \mathbb{E}\left[\left(\sup_{u\leq T}\left(t-g_{\mu}\left(t\right)\right)\right)^{k\mu}\right] &\leq&  \dfrac{2-k}{\left(1-k\right)2^{k\mu}\Gamma\left(1+\mu\right)^{k}}\mathbb{E}\left[L_{T}^{k}\right]\\
   &\leq& \left(\dfrac{2-k}{1-k}\right)^{2}\mathbb{E}\left[\left(\sup_{u\leq T}\left(t-g_{\mu}\left(t\right)\right)\right)^{k\mu}\right].
\end{eqnarray*}There are also other types of inequalities for local times of
Az\'{e}ma martingales ($\mu=\dfrac{1}{2}$) obtained by Chao and Chou
\cite{chaochou}.
\end{rem}
Now, we shall give some inequalities, obtained by Yor
(\cite{yorinegal}), for local submartingales of the class
$(\Sigma)$. In \cite{yorinegal}, the inequalities are given for
continuous local submartingales; the proof in fact applies to local
submartingales of the class $(\Sigma)$ with only negative jumps if
one uses the version of Skorokhod's reflection lemma for functions
with jumps, as explained in \cite{Mckean} or \cite{ashyordoob}.

First, for $X=N+A$ of the class $(\Sigma)$, and $k\in(0,\infty)$, we
define:
$$\left\|X\right\|_{H^{k}}=\left\|[N,N]_{\infty}^{1/2}\right\|_{L^{k}}+\left\|A_{\infty}\right\|_{L^{k}}.$$For
$k\geq1$, it defines a norm. Now, we state a slight generalization
of Yor's inequalities (in this version negative jumps are allowed):
\begin{prop}[Yor \cite{yorinegal}]
Let $X$ be a local submartingale of class $(\Sigma)$, which has only
negative jumps. Then the following holds:
\begin{enumerate}
\item for all $k\in[1,\infty)$,
$\left\|A_{\infty}\right\|_{L^{k}}\leq\left\|N_{\infty}^{*}\right\|_{L^{k}}$,
and the quantities $\left\|X\right\|_{H^{k}}$,
$\left\|X_{\infty}^{*}\right\|_{L^{k}}$ and
$\left\|N\right\|_{H^{k}}$ are equivalent.
\item if $N$ is continuous, then the above statement holds for all $k\in(0,\infty)$, and if $k\in(0,1)$, then $\left\|A_{\infty}\right\|_{L^{k}}$ is
equivalent to the three former quantities.
\end{enumerate}
\end{prop}
\section{The distribution of $A_{\infty}$}
In this section, we shall compute the law of $A_{\infty}$, for a
large class of processes in $(\Sigma D)$, and then we shall see that
we are able to recover some well known results (for example,
formulas for stopped diffusions \cite{lehoczky}).
\begin{thm}\label{loideA}
Let $X$ be a process of class $(\Sigma D)$, and define
$$\lambda\left(x\right)\equiv\mathbb{E}\left[X_{\infty}|A_{\infty}=x\right].$$
Assume that $\lambda\left(A_{\infty}\right)\neq0.$ Then, if we note
$$b\equiv\inf\left\{u:\;\mathbb{P}\left(A_{\infty}\geq u\right)=0\right\},$$we
have:
\begin{equation}\label{loiAinfini}
\mathbb{P}\left(A_{\infty}>x\right)=\exp\left(-\int_{0}^{x}\dfrac{dz}{\lambda\left(z\right)}\right);\;x<b.
\end{equation}
\end{thm}
\begin{rem}
The techniques we shall use for the proof of Theorem \ref{loideA}
are very close to those already used by Az\'{e}ma and Yor
(\cite{AY}) and more recently by Vallois (\cite{vallois}) in the
study of the law of the maximum of a continuous and uniformly
integrable martingale. Our results here allow for some
discontinuities (look at the examples in Section 2).
\end{rem}
\begin{proof}
Let $f$ be a bounded Borel function with compact support; from
Corollary \ref{corolfbornee},
$$\mathbb{E}\left[F\left(A_{\infty}\right)\right]=\mathbb{E}\left[X_{\infty}f\left(A_{\infty}\right)\right],$$where
$F\left(x\right)=\int_{0}^{x}dzf\left(z\right)$. Now, conditioning
the right hand side with respect to $A_{\infty}$ yields:
\begin{equation}\label{a}
    \mathbb{E}\left[F\left(A_{\infty}\right)\right]=\mathbb{E}\left[\lambda\left(A_{\infty}\right)f\left(A_{\infty}\right)\right].
\end{equation}
Now, since $\lambda\left(A_{\infty}\right)>0$, if
$\nu\left(dx\right)$ denotes the law of $A_{\infty}$, and
$\overline{\nu}\left(x\right)\equiv\nu\left(\left[x,\infty\right)\right)$,
(\ref{a}) implies:
\begin{equation*}
  \int_{0}^{\infty}dzf\left(z\right)\overline{\nu}\left(z\right)=
  \int_{0}^{\infty}\nu\left(dz\right)f\left(z\right)\lambda\left(z\right),
\end{equation*}and consequently,
\begin{equation}\label{b}
    \overline{\nu}\left(z\right)dz=\lambda\left(z\right)\nu\left(dz\right).
\end{equation}
Recall that $b\equiv\inf\left\{u:\;\mathbb{P}\left(A_{\infty}\geq
u\right)=0\right\}$; hence for $x<b$,
$$\int_{0}^{x}\dfrac{dz}{\lambda\left(z\right)}=\int_{0}^{x}\dfrac{\nu\left(dz\right)}{\overline{\nu}\left(z\right)}\leq\dfrac{1}{\overline{\nu}\left(z\right)}<\infty,$$
and integrating (\ref{b}) between $0$ and $x$, for $x<b$ yields:
$$\overline{\nu}\left(x\right)=\exp\left(-\int_{0}^{x}\dfrac{dz}{\lambda\left(z\right)}\right),$$and
the result of the Theorem follows easily.
\end{proof}
\begin{rem}
We note that:
$$b=\inf\left\{x:\;\int_{0}^{x}\dfrac{dz}{\lambda\left(z\right)}=+\infty\right\},$$with
the usual convention that $\inf\emptyset=\infty$.
\end{rem}
\begin{rem}
We took the hypothesis $\lambda\left(A_{\infty}\right)\neq0$ to
avoid technicalities. The result stated in Theorem \ref{loideA} is
general enough for our purpose; the reader interested in the general
case can refer to \cite{vallois} where the special case
$X_{t}=S_{t}-M_{t}$, with $M$ continuous, is dealt with in depth.
\end{rem}
Now, we give a series of interesting corollaries.
\begin{cor}
Let $\left(M_{t}\right)$ be a uniformly integrable martingale, with
only negative jumps, and such that
$\mathbb{E}\left(S_{\infty}\right)<\infty$, and
$M_{\infty}<S_{\infty}$. Then,
$$\mathbb{P}\left(S_{\infty}>x\right)=\exp\left(-\int_{0}^{x}\dfrac{dz}{z-\alpha\left(z\right)}\right),$$where
$$\alpha\left(x\right)=\mathbb{E}\left[M_{\infty}|S_{\infty}=x\right].$$
\end{cor}
\begin{rem}
This result, when the martingale $M$ is continuous, is a special
case of a more general result by Vallois \cite{vallois} where
$S_{\infty}-\alpha\left(S_{\infty}\right)$ can vanish. The reader
can also refer to \cite{rogers} for a discussion on the law of the
maximum of a martingale and its terminal value.
\end{rem}
\begin{cor}\label{loiexppourAinfini}
Let $X$ be a process of class $(\Sigma D)$, such that:
$$\lim_{t\rightarrow\infty}X_{t}=a>0,\;a.s.$$Then,
$$\mathbb{P}\left(A_{\infty}>x\right)=\exp\left(-\dfrac{x}{a}\right),$$i.e. $A_{\infty}$
is distributed as a random variable with exponential law of
parameter $\frac{1}{a}$.
\end{cor}
\begin{proof}
It is a consequence of Theorem \ref{loideA} with
$\lambda\left(x\right)\equiv a$.
\end{proof}
Let $\left(M_{t}\right)$ be a continuous local martingale such that
$<M>_{\infty}=\infty,\;\mathrm{a.s.}$; let $T_{1}=\inf\left\{t\geq
0:\;M_{t}=1\right\}$. Then an application of Tanaka's formula and
Corollary \ref{loiexppourAinfini} shows that
$\dfrac{1}{2}L_{T_{1}}\left(M\right)$ follows the standard
exponential law ($L_{T_{1}}\left(M\right)$ denotes the local time at
$0$ of the local martingale $M$). The result also applies to Bessel
processes of dimension $2(1-\mu)$, with $\mu\in(0,1)$: taking
$T_{1}=\inf\left\{t\geq 0:\;R_{t}=1\right\}$, we have that
$L_{T_{1}}\left(R\right)$ follows the standard exponential law. The
result also applies to $\left(t-g_{\mu}\left(t\right)\right)$.
\begin{cor}
Let $\left(M_{t}\right)$ be a continuous martingale such that
$\lim_{t\rightarrow\infty}M_{t}=M_{\infty}$ exists and
$|M_{\infty}|>0$. Then, if $\left(L_{t}\right)$ denotes its local
time at $0$, we have:
$$\mathbb{P}\left(L_{\infty}>x\right)=\exp\left(-\int_{0}^{x}\dfrac{dz}{\mathbb{E}\left[|M_{\infty}|\;|L_{\infty}=z\right]}\right).$$
\end{cor}
\begin{cor}\label{AarreteenT}
Let $X$ be of the class $(\Sigma)$ with only negative jumps and
$A_{\infty}=\infty$ and let $\varphi$ be a nonnegative locally
bounded Borel function such that
$\int_{0}^{\infty}dx\varphi\left(x\right)=\infty$. Define the
stopping time $T$ as:
$$T\equiv\inf\left\{t:\;\varphi\left(A_{t}\right)X_{t}=1\right\}.$$Then
$T<\infty,\;a.s.$ and
$$\mathbb{P}\left(A_{T}>x\right)=\exp\left(-\int_{0}^{x}dz\varphi\left(z\right)\right).$$
\end{cor}
\begin{proof}
The fact that $T<\infty,\;a.s.$ is a consequence of Theorem
\ref{loideA} and the rest follows from Theorem \ref{loideA} with
$\lambda\left(x\right)=\frac{1}{\varphi\left(x\right)}$.
\end{proof}
The following variant of Corollary \ref{AarreteenT} is sometimes
useful:
\begin{cor}\label{AarreteenTbis}
Let $X$ be of the class $(\Sigma)$ with only negative jumps and
$A_{\infty}=\infty$. Let $\psi$ be a nonnegative Borel function such
that $\dfrac{1}{\psi}$ is locally bounded and
$\int_{0}^{\infty}\frac{dx}{\psi\left(x\right)}=\infty$. Define the
stopping time $T$ as:
$$T\equiv\inf\left\{t:\;X_{t}\geq\psi\left(A_{t}\right)\right\}.$$Then
$T<\infty,\;a.s.$ and
$$\mathbb{P}\left(A_{T}>x\right)=\exp\left(-\int_{0}^{x}\dfrac{dz}{\psi\left(z\right)}\right).$$
\end{cor}

Now, we shall apply the previous results to compute the maximum of a
stopped continuous diffusion process; in particular, we are able to
recover a formula discovered first by Lehoczky (\cite{lehoczky}).
More precisely, let $\left(Y_{t}\right)$ be a continuous diffusion
process, with $Y_{0}=0$. Let us assume further that $Y$ is
recurrent; then, from the general theory of diffusion processes,
there exists a unique continuous and strictly increasing function
$s$, with $s\left(0\right)=0$, $\lim_{x\rightarrow
+\infty}s\left(x\right)=+\infty$, $\lim_{x\rightarrow
-\infty}s\left(x\right)=-\infty$, such that $s\left(Y_{t}\right)$ is
a continuous local martingale. Let $\theta$ be a Borel function with
$\theta\left(x\right)>0,\;\forall x>0$. Define the stopping time
$$T\equiv\inf\left\{t:\;\overline{Y}_{t}-Y_{t}\geq\theta\left(\overline{Y}_{t}\right)\right\},$$where
$\overline{Y}_{t}=\sup_{u\leq t}Y_{u}$.
\begin{prop}\label{thlehocgen}
Let us assume that $T<\infty,\;a.s.$. Then, the law of
$\overline{Y}_{T}$ is given by:
\begin{equation}\label{Lehocgen}
    \mathbb{P}\left(\overline{Y}_{T}>x\right)=\exp\left(-\int_{0}^{x}\dfrac{ds\left(z\right)}{s\left(z\right)-s\left(z-\theta\left(z\right)\right)}\right).
\end{equation}In particular, when $\theta\left(x\right)\equiv a$,
with $a>0$, we have:
$$\mathbb{P}\left(\overline{Y}_{T}>x\right)=\exp\left(-\int_{0}^{x}\dfrac{ds\left(z\right)}{s\left(z\right)-s\left(z-a\right)}\right).$$
\end{prop}
\begin{proof}
First, we note that
$$\overline{Y}_{T}-Y_{T}=\theta\left(\overline{Y}_{T}\right).$$
Indeed, considering the process
$K_{t}=\dfrac{\overline{Y}_{t}-Y_{t}}{\theta\left(\overline{Y}_{t}\right)}$,
we have $T=\inf\left\{t:\;K_{t}=1\right\}$, and each time
$\theta\left(\overline{Y}_{t}\right)$ jumps corresponds to an
increase time for $\overline{Y}_{t}$, and since at such a time
$\overline{Y}_{t}-Y_{t}=0$, $K$ is in fact continuous and $K_{T}=1$.

Now, let us define $X$ of the class $(\Sigma)$ by:
$$X_{t}=s\left(\overline{Y}_{t}\right)-s\left(Y_{t}\right).$$ From
the remark above, we have:
$$X_{T}=s\left(\overline{Y}_{T}\right)-s\left(\overline{Y}_{T}-\theta\left(\overline{Y}_{T}\right)\right).$$Now,
with the notations of Theorem \ref{loideA}, we have:
$$\lambda\left(x\right)=\mathbb{E}\left[X_{T}|\overline{Y}_{T}=x\right]=x-s\left(s^{-1}\left(x\right)-\theta\left(s^{-1}\left(x\right)\right)\right),$$
and consequently, since $\lambda\left(x\right)>0$, from Theorem
\ref{loideA}, we have:
$$\mathbb{P}\left(s\left(\overline{Y}_{T}\right)>x\right)=\exp\left(-\int_{0}^{x}\dfrac{dz}{z-s\left(s^{-1}\left(z\right)-\theta\left(s^{-1}\left(z\right)\right)\right)}\right),$$and
thus
$$\mathbb{P}\left(\overline{Y}_{T}>x\right)=\exp\left(-\int_{0}^{s\left(x\right)}\dfrac{dz}{z-s\left(s^{-1}\left(z\right)-\theta\left(s^{-1}\left(z\right)\right)\right)}\right).$$
Now, making the change of variable $z=s\left(u\right)$ gives the
desired result.
\end{proof}
Now, if $Y$ is of the form:
$$dY_{t}=b\left(Y_{t}\right)dt+\sigma\left(Y_{t}\right)dB_{t},$$where $B$
is a standard Brownian Motion, and the coefficients $b$ and $\sigma$
chosen such that uniqueness, existence and recurrence hold, then
$$s\left(x\right)=\int_{0}^{x}dy\exp\left(-\beta\left(y\right)\right),$$with
$$\beta\left(x\right)=2\int_{0}^{x}dy\dfrac{b\left(y\right)}{\sigma^{2}\left(y\right)}.$$Under
these assumptions, (\ref{Lehocgen}) takes the following form:
\begin{cor}[Lehoczky \cite{lehoczky}]
With the assumptions and notations of Proposition \ref{thlehocgen},
and the notations above, we have:
\begin{equation}
    \mathbb{P}\left(\overline{Y}_{T}>x\right)=\exp\left(-\int_{0}^{x}\dfrac{\exp\left(-\beta\left(z\right)\right)dz}{\int_{z-\theta\left(z\right)}^{z}\exp\left(-\beta\left(u\right)\right)du}\right).
\end{equation}In the special case when $\theta\left(x\right)\equiv a$,
with $a>0$, we have:
\begin{equation}
    \mathbb{P}\left(\overline{Y}_{T}>x\right)=\exp\left(-\int_{0}^{x}\dfrac{\beta\left(z\right)dz}{\int_{z-a}^{z}\beta\left(u\right)du}\right).
\end{equation}
\end{cor}\bigskip
To conclude this section, we mention a class of stochastic
processes, which are solution of a stochastic differential equation
of Skorokhod type, with reflecting boundary condition at $0$, and
which look very similar to the local submartingales of class
$(\Sigma)$. These processes play an important role in the extension
of Pitman's theorem to one dimensional diffusion processes in Saisho
and Tanemura's work \cite{saichotanemura}; they also appeared
earlier in the works of Chaleyat-Maurel and El Karoui (see their
paper in \cite{asterisque}).

More precisely, let
$\left(\Omega,\mathcal{F},\left(\mathcal{F}_{t}\right),\mathbb{P}\right)$
be a filtered probability space and $\left(B_{t}\right)$ an
$\left(\mathcal{F}_{t}\right)$ Brownian Motion. Let
$\sigma,\;b:\mathbb{R}\rightarrow\mathbb{R}$ be Lipschitz continuous
functions, and assume that $\sigma\left(x\right)>0,\;\forall
x\in\mathbb{R}$. Consider the stochastic differential equation of
Skorokhod type:
\begin{equation}\label{diffusionreflechie}
Y_{t}=\int_{0}^{t}\sigma\left(Y_{u}+L_{u}\right)dB_{u}+\int_{0}^{t}b\left(Y_{u}+L_{u}\right)du+L_{t},
\end{equation}where $Y$ and $L$ should be found under the conditions:
\begin{itemize}
\item $Y$ is $\left(\mathcal{F}_{t}\right)$, continuous and
$Y_{t}\geq0$;
\item $L$ is continuous, nondecreasing, $L_{0}=0$ and
$L_{t}=\int_{0}^{t}\mathbf{1}_{\left\{0\right\}}\left(Y_{u}\right)dL_{u}$.
\end{itemize}It is proved in \cite{saichotanemura} that this
equation has a unique solution. Although $Y$ is not a process of the
class $(\Sigma)$, we shall apply the previous methods to compute the
law of $L$ conveniently stopped.

We let as before $s$ denote the function:
$$s\left(x\right)=\int_{0}^{x}dy\exp\left(-\beta\left(y\right)\right),$$with
$$\beta\left(x\right)=2\int_{0}^{x}dy\dfrac{b\left(y\right)}{\sigma^{2}\left(y\right)}.$$
If $\left(Y_{t},L_{t}\right)$ is the solution of
(\ref{diffusionreflechie}), define:
\begin{eqnarray*}
  X_{t} &=& s\left(Y_{t}+L_{t}\right)-s\left(L_{t}\right), \\
  N_{t} &=& \int_{0}^{t}\sigma s'\left(Y_{u}+L_{u}\right)dB_{u}, \\
  A_{t} &=& s\left(L_{t}\right).
\end{eqnarray*}
We shall need the following lemma:
\begin{lem}
The following equality holds:
$$X_{t}=N_{t}+L_{t},$$and $X$ is of the class $\left(\Sigma\right)$.
\end{lem}
\begin{proof}
An application of It\^{o}'s formula yields:
\begin{equation}\label{zaza}
    s\left(Y_{t}+L_{t}\right)=\int_{0}^{t}\sigma
s'\left(Y_{u}+L_{u}\right)dB_{u}+2\int_{0}^{t}s'\left(Y_{u}+L_{u}\right)dL_{u}.
\end{equation}
Now, since  $\left(dL_{t}\right)$ is carried by the set of zeros of
$Y$, we have:
$$\int_{0}^{t}s'\left(Y_{u}+L_{u}\right)dL_{u}=\int_{0}^{t}s'\left(L_{u}\right)dL_{u}=s\left(L_{t}\right)=A_{t},$$and
consequently, from (\ref{zaza}), we have:
$$s\left(Y_{t}+L_{t}\right)-s\left(L_{t}\right)=\int_{0}^{t}\sigma
s'\left(Y_{u}+L_{u}\right)dB_{u}+s\left(L_{t}\right),$$that is:
$$X_{t}=N_{t}+L_{t}.$$It follows easily from the fact that $s$ is
continuous and strictly increasing that $A$ is increasing and that
$\left(dA_{t}\right)$ is carried by the set of zeros of $X$. Hence
$X$ is of the class $\left(\Sigma\right)$.
\end{proof}
Now, we can state  an  analogue of Lehoczky's result for the pair:
$\left(Y_{t},L_{t}\right)$:
\begin{prop}
Let $T$ be the stopping time
$$T\equiv\inf\left\{t:\;Y_{t}\geq\theta\left(L_{t}\right)\right\},$$where
$\theta:\mathbb{R}_{+}\rightarrow\mathbb{R}_{+}$ is a Borel function
such that $\theta\left(x\right)>0,\;\forall x\geq 0$. Then, if
$T<\infty,\;a.s.$,
$$\mathbb{P}\left(L_{T}>x\right)=\exp\left(-\int_{0}^{x}\dfrac{\exp\left(-\beta\left(z\right)\right)dz}{\int_{z}^{z+\theta\left(z\right)}\exp\left(-\beta\left(u\right)\right)du}\right).$$
In the special case $\theta\left(x\right)\equiv a$, for some $a>0$,
we have:
$$\mathbb{P}\left(L_{T}>x\right)=\exp\left(-\int_{0}^{x}\dfrac{\exp\left(-\beta\left(z\right)\right)dz}{\int_{z}^{z+a}\exp\left(-\beta\left(u\right)\right)du}\right).$$
\end{prop}
\begin{proof}The proof follows exactly the same line as the proof of
Theorem \ref{thlehocgen}, so we just give the main steps.  Here
again, we have:
$$Y_{T}=\theta\left(L_{T}\right),$$and consequently,
$$X_{T}=s\left(\theta\left(L_{T}\right)+L_{T}\right)-s\left(L_{T}\right),$$and
hence:
$$\lambda\left(x\right)=s\left(s^{-1}\left(x\right)+\theta\left(s^{-1}\left(x\right)\right)\right)-x.$$The
end of the proof is now exactly the same as that of Theorem
\ref{thlehocgen}.
\end{proof}
\begin{rem}
If the ratio
$\dfrac{b\left(x\right)}{\sigma^{2}\left(x\right)}\equiv\gamma$ is
constant, then  $L_{T}$, with $T\equiv\inf\left\{t:\;Y_{t}\geq
a\right\}$, is exponentially distributed:
$$\mathbb{P}\left(L_{T}>x\right)=\exp\left(-\frac{2\gamma x}{1-\exp\left(-2\gamma a\right)}\right).$$
\end{rem}
\section{The Skorokhod embedding problem for non atomic probability
measures on $\mathbb{R}_{+}$ } The previous results can now be used
to solve the Skorokhod stopping problem for non atomic probability
measures on $\mathbb{R}_{+}$. The literature on the Skorokhod
stopping problem is vast (see \cite{jansurvey}) and the aim here is
not to give very specialized results on this topic, but rather
illustrate a general methodology which allows to deal with a wide
variety of stochastic processes, which may even be discontinuous. If
the reader is interested in solving the problem for measures which
have atoms, then he can refer to the recent paper of Ob\l\'{o}j and
Yor \cite{jan}, which inspired the ideas in the sequel, and which
explains how to deal with such probability measures. Quite
remarkably, the stopping times they propose for the reflected
Brownian Motion and the age process of Brownian excursions (hence
solving explicitly for the first time the Skorokhod embedding
problem for a discontinuous process) can be used to solve the
Skorokhod stopping problem for local submartingales of the class
$(\Sigma)$, with only negative jumps and with $A_{\infty}=\infty$.
In particular, we will be able to solve the Skorokhod stopping
problem for any Bessel process of dimension $2(1-\mu)$, with
$\mu\in(0,1)$, and its corresponding age process. In \cite{jan}, the
authors mention that they have presented detailed arguments using
excursion theory, but parallel arguments, using martingale theory,
are possible as well. As far as we are concerned here, we shall
detail arguments based on martingale theory, to show that the
stopping times they propose have a universal feature (independent of
any Markov or scaling property).\bigskip

More precisely, let $\vartheta$ be a probability measure on
$\mathbb{R}_{+}$, which has no atoms, and let $X$ be a local
submartingale of the class $(\Sigma)$, with only negative jumps and
such that $\lim_{t\rightarrow\infty}A_{t}=\infty$. Our aim is to
find a stopping time $T_{\vartheta}$, such that the law of
$X_{T_{\vartheta}}$ is $\vartheta$, and which coincides with the
stopping time proposed by Obl\'{o}j and Yor \cite{jan} when
$X_{t}=B_{t}$ or $X_{t}=\sqrt{t-g_{1/2}\left(t\right)}$.

In the sequel, we write
$\overline{\vartheta}\left(x\right)=\vartheta\left(\left[x,\infty\right)\right)$
for the tail of $\vartheta$, and
$a_{\vartheta}=\sup\left\{x\geq0:\;\overline{\vartheta}\left(x\right)=1\right\}$,
and
$b_{\vartheta}=\inf\left\{x\geq0:\;\overline{\vartheta}\left(x\right)=0\right\}$,
$-\infty\leq a_{\vartheta}\leq b_{\vartheta}\leq\infty$,
respectively, for the lower and upper bound of the support of
$\vartheta$. Now, following \cite{jan}, we introduce the dual
Hardy-Littlewood function
$\psi_{\vartheta}:[0,\infty)\rightarrow[0,\infty)$ through:
$$\psi_{\vartheta}\left(x\right)=\int_{\left[0,x\right]}\dfrac{z}{\overline{\vartheta}\left(z\right)}d\vartheta\left(z\right),\;a_{\vartheta}\leq x<b_{\vartheta},$$
and $\psi_{\vartheta}\left(x\right)=0$ for $0\leq x<a_{\vartheta}$,
and $\psi_{\vartheta}\left(x\right)=\infty$ for $x\geq
b_{\vartheta}$. The function $\psi_{\vartheta}$ is continuous and
increasing, and we can define its right continuous inverse
$$\varphi_{\vartheta}\left(z\right)=\inf\left\{x\geq0:\;\psi_{\vartheta}\left(x\right)>z\right\},$$which
 is strictly increasing. Now, we can state
the main result of this subsection:
\begin{thm}\label{skorod}
Let $X$ be a local  submartingale of the class $(\Sigma)$, with only
negative jumps and such that $A_{\infty}=\infty$. The stopping time
\begin{equation}\label{tpsarretsko}
    T_{\vartheta}=\inf\left\{x\geq0:\;X_{t}\geq\varphi_{\vartheta}\left(A_{t}\right)\right\}
\end{equation}is $a.s.$ finite and solves the Skorokhod embedding
problem for $X$, i.e. the law of $X_{T_{\vartheta}}$ is $\vartheta$.
\end{thm}
\begin{proof}
Part of the arguments that follow are tailored on those of
Ob\l\'{o}j and Yor \cite{jan}, but at some stage, they use time
change techniques, which would not deal with the case of negative
jumps.

First, we note that (see \cite{jan}):
\begin{eqnarray*}
  \int_{0}^{x}\dfrac{dz}{\varphi_{\vartheta}\left(z\right)} &<& \infty,\;\mathrm{for\ }0\leq x< b_{\vartheta} \\
  \int_{0}^{\infty}\dfrac{dz}{\varphi_{\vartheta}\left(z\right)} &=&
  \infty,
\end{eqnarray*}and consequently, from Theorem \ref{estimationavecA} and Remark \ref{phiategalx},
$T_{\vartheta}<\infty,\;a.s.$ and
$X_{T_{\vartheta}}=\varphi_{\vartheta}\left(A_{T_{\vartheta}}\right)$.

Now, let $h:\mathbb{R}_{+}\rightarrow\mathbb{R}_{+}$ be a strictly
decreasing function, locally bounded, and such that
$\int_{0}^{\infty}dzh\left(z\right)=\infty$. From Theorem
\ref{caracteriation1}, $h\left(A_{t}\right)X_{t}$ is again of the
class $(\Sigma)$, and its increasing process is
$\int_{0}^{A_{t}}dzh\left(z\right)\equiv H\left(A_{t}\right)$. Now,
define the stopping time:
$$R_{h}=\inf\left\{t:\;h\left(A_{t}\right)X_{t}=1\right\},$$which is
finite almost surely from Theorem \ref{estimationavecA}. Since
$\lim_{t\rightarrow\infty}H\left(A_{t}\right)=\infty$, from
Corollary \ref{loiexppourAinfini}, $H\left(A_{R_{h}}\right)$ is
distributed as a random variable $\mathbf{e}$ which follows the
standard exponential law, and hence: $$A_{R_{h}}\
\stackrel{\mbox{\small (law)}}{=}\ \
H^{-1}\left(\mathbf{e}\right),$$and consequently:
\begin{equation}\label{loidexrh}
    X_{R_{h}}\
\stackrel{\mbox{\small (law)}}{=}\ \
\dfrac{1}{h\left(H^{-1}\left(\mathbf{e}\right)\right)}.
\end{equation}

Now, we investigate the converse problem, that is given a
probability measure $\vartheta$, we want to find $h$ such that
$X_{R_{h}}\ \stackrel{\mbox{\small (law)}}{=}\ \ \vartheta$. From
(\ref{loidexrh}), we deduce (recall that $\vartheta$ has no atoms):
\begin{eqnarray}
  \overline{\vartheta}\left(x\right) &=& \mathbb{P}\left(\dfrac{1}{h\left(H^{-1}\left(\mathbf{e}\right)\right)}>x\right)=\mathbb{P}\left(h\left(H^{-1}\left(\mathbf{e}\right)\right)<\dfrac{1}{x}\right) \notag\\
   &=& \mathbb{P}\left(H^{-1}\left(\mathbf{e}\right)>h^{-1}\left(\dfrac{1}{x}\right)\right) \notag\\
   &=&
   \mathbb{P}\left(\mathbf{e}>H\left(h^{-1}\left(\dfrac{1}{x}\right)\right)\right)=\exp\left(-H\left(h^{-1}\left(\dfrac{1}{x}\right)\right)\right).
   \label{sko1}
\end{eqnarray}Now, differentiating the last equality yields:
\begin{eqnarray*}
  -d\overline{\vartheta}\left(x\right) &=& \overline{\vartheta}\left(x\right)\left[h\left(h^{-1}\left(\dfrac{1}{x}\right)\right)\right]d\left(h^{-1}\left(\dfrac{1}{x}\right)\right) \\
   &=&
   \dfrac{\overline{\vartheta}\left(x\right)}{x}d\left(h^{-1}\left(\dfrac{1}{x}\right)\right),
\end{eqnarray*}hence
\begin{equation}\label{sko2}
    d\left(h^{-1}\left(\dfrac{1}{x}\right)\right)=x\dfrac{d\vartheta\left(x\right)}{\overline{\vartheta}\left(x\right)}.
\end{equation}
Consequently, we have:
$$h^{-1}\left(\dfrac{1}{x}\right)=\int_{0}^{x}\dfrac{z}{\overline{\vartheta}\left(z\right)}d\vartheta\left(z\right)=\psi_{\vartheta}\left(x\right),$$
and
$$R_{h}=\inf\left\{t:\;X_{t}\geq\ \varphi_{\vartheta}\left(A_{t}\right)\right\}=T_{\varphi},$$and
the proof of the theorem follows easily.
\end{proof}
\begin{rem}
One can easily check that
$$X_{t\wedge
T\vartheta}\;\mathrm{is\;uniformly\;integrable}\;\Leftrightarrow
\int_{0}^{\infty}xd\vartheta\left(x\right)<\infty.$$
\end{rem}
Theorem \ref{skorod} shows that the stopping times proposed by
Ob\l\'{o}j and Yor for the pair $\left(|B_{t}|,\ell_{t}\right)$ and
$\left(\sqrt{\frac{\pi}{2}\left(t-g_{1/2}\left(t\right)\right)},\ell_{t}\right)$
have a universal aspect in the sense that they apply to a very wide
class of processes. For example, Theorem \ref{skorod} provides us
with a solution for the Skorokhod stopping problem for Bessel
processes of dimension $2(1-\mu)$, with $\mu\in(0,1)$, whilst these
processes are not semimartingales for $\mu>1/2$. With Theorem
\ref{skorod}, we can also solve the Skorokhod stopping problem for
the age processes associated with those Bessel processes. In fact,
the Skorokhod embedding problem is solved for powers of these
processes, but one can then easily deduce the stopping time for the
process itself. More precisely:
\begin{cor}
Let $\left(R_{t}\right)$ be a Bessel process, starting from $0$, of
dimension $2(1-\mu)$, with $\mu\in(0,1)$, and $\left(L_{t}\right)$
its local time. Let $\vartheta$ be a non atomic probability measure
on $\mathbb{R}_{+}$; then the stopping time
$$T_{\vartheta}=\inf\left\{x\geq0:\;R_{t}\geq\varphi_{\vartheta}^{1/2\mu}\left(L_{t}\right)\right\}$$
is $a.s.$ finished and solves the Skorokhod embedding problem for
$R_{t}^{2\mu}$. In particular, when $\mu=1/2$, we obtain the
stopping time proposed by Obl\'{o}j and Yor for $|B_{t}|$.
\end{cor}
\begin{cor}
With the notations and assumptions of the previous corollary, let
$$g_{\mu}\left(t\right)\equiv\sup\left\{u\leq t:\;R_{u}=0\right\}.$$
Recall that in the filtration
$\left(\mathcal{G}_{t}\equiv\mathcal{F}_{g_{\mu}\left(t\right)}\right)$
of the zeros of the Bessel process $R$,
$\left(t-g_{\mu}\left(t\right)\right)^{\mu}$, is a submartingale of
class $(\Sigma)$ whose increasing process in its Doob-Meyer
decomposition is given by: $A_{t}\equiv
\dfrac{1}{2^{\mu}\Gamma\left(1+\mu\right)}L_{t}$. Consequently, the
stopping time
$$T_{\vartheta}=\inf\left\{x\geq0:\;t-g_{\mu}\left(t\right)\geq\varphi_{\vartheta}^{1/\mu}\left(\dfrac{1}{2^{\mu}\Gamma\left(1+\mu\right)}L_{t}\right)\right\}$$
is $a.s.$ finished and solves the Skorokhod embedding problem for
$\left(t-g_{\mu}\left(t\right)\right)^{\mu}$. In particular, when
$\mu=1/2$, we obtain the stopping time proposed by Ob\l\'{o}j and
Yor for $\sqrt{t-g_{1/2}\left(t\right)}$, the age process of a
standard Brownian Motion.
\end{cor}
\begin{rem}
Our methodology does not allow us to solve the Skorokhod stopping
problem for local submartingales of the class $(\Sigma)$ which have
positive jumps. It would thus be interesting to discover a method
which would lead us to a solution for these processes.
\end{rem}\bigskip

One can also use the ideas contained in the proof of Theorem
\ref{skorod} to give another family of stopping times which solve
the Skorokhod stopping problem. The stopping time we shall propose
was obtained first by Obl\'{o}j and Yor for the absolute value of
the standard Brownian Motion. Here again, we shall prove that this
stopping time solves the Skorokhod embedding problem for a much
wider class of stochastic processes. Indeed, if in the proof of
Theorem \ref{skorod} we had supposed $h$ strictly increasing instead
of decreasing, then (\ref{sko1}) would have equaled
$\vartheta\left(x\right)\equiv\vartheta\left(\left[0,x\right]\right)$
instead of $\overline{\vartheta}\left(x\right)$ and (\ref{sko2})
would read:
$$d\left(h^{-1}\left(\dfrac{1}{x}\right)\right)=-x\dfrac{d\vartheta\left(x\right)}{\vartheta\left(x\right)}.$$Therefore,
we have the following proposition:
\begin{prop}
Let $X$ be a local  submartingale of the class $(\Sigma)$, with only
negative jumps and with $A_{\infty}=\infty$. Let $\vartheta$ be a
probability measure on $\mathbb{R}_{+}$, without atoms, and define
$$\overline{\psi}_{\vartheta}\left(x\right)=\int_{\left(x,\infty\right)}\dfrac{xd\vartheta\left(x\right)}{\vartheta\left(x\right)},$$and
$$\overline{\varphi}_{\vartheta}\left(x\right)=\inf\left\{z\geq0:\;\overline{\psi}_{\vartheta}\left(z\right)>x\right\},\;x<b_{\vartheta},$$and
$\overline{\varphi}_{\vartheta}\left(x\right)=0,\;\forall x\geq
b_{\vartheta}$. Then the stopping time
$$\overline{T}_{\vartheta}=\inf\left\{t\geq0:\;X_{t}\geq\overline{\varphi}_{\vartheta}\left(A_{t}\right)\right\}$$embeds
$\vartheta$, i.e.
$$X_{\overline{T}_{\vartheta}}\
\stackrel{\mbox{\small (law)}}{=}\ \ \vartheta.$$
\end{prop}
\begin{rem}
Here again, the stopping time $\overline{T}_{\vartheta}$ embeds
$\vartheta$ for the Bessel processes of dimension $2(1-\mu)$,
$\mu\in(0,1)$, and the corresponding age processes. The fact that
$\overline{T}_{\vartheta}$ works for the Brownian age process was
not noticed by Obl\'{o}j and Yor.
\end{rem}\bigskip

Let us conclude this subsection with two simple examples.
\begin{ex}[Exponential law]
Let $\vartheta\left(dx\right)=\rho\exp\left(-\rho
x\right)\mathbf{1}_{\mathbb{R}_{+}}\left(x\right)dx,\;\rho>0$. We
then have:
$$\varphi_{\vartheta}\left(x\right)=\sqrt{\dfrac{2}{\rho}x}.$$The
stopping time takes the following form:
$$T_{\vartheta}=\inf\left\{t:\;X_{t}\geq\sqrt{\dfrac{2}{\rho}A_{t}}\right\}.$$
\end{ex}
\begin{ex}[Uniform law]
Let
$\vartheta\left(dx\right)=\dfrac{1}{b}\mathbf{1}_{\left[0,b\right]}\left(x\right)dx$.
Some simple calculations give:
$$\overline{\psi}_{\vartheta}\left(x\right)=\left(b-x\right),$$and
therefore
$$\overline{T_{\vartheta}}=\inf\left\{t:\;X_{t}+A_{t}\geq b\right\}.$$ The reader can easily check that the corresponding formula for
$T_{\vartheta}$ is not so simple.
\end{ex}
\section*{Acknowledgements}
I am very grateful to my supervisor Marc Yor for many helpful
discussions and for correcting earlier versions of this paper.

\end{document}